\documentclass[12pt]{article}
\usepackage{latexsym}
\usepackage{amsfonts}
\usepackage{amsmath}
\usepackage{amssymb}
\usepackage{graphicx}
\usepackage{latexsym}
\usepackage{amsfonts}
\usepackage{graphicx}
\usepackage{psfrag}

\input{tcilatex}
\begin{document}

\qquad %\setcounter{section}{-1}

\thispagestyle{empty}

\begin{center}
{\Large \textbf{\ NOTE ON BOUNDS FOR EIGENVALUES USING TRACES }}

\bigskip

R. Sharma, R. Kumar and R. Saini

Department of Mathematics, Himachal Pradesh University, Shimla 171 005, India%
\\[0pt]
email: rajesh\underline{~~}hpu\underline{~~}math@yahoo.co.in

\bigskip
\end{center}

\vskip1.5in \noindent \textbf{Abstract. }We show that various old and new
bounds involving eigenvalues of a complex $n\times n$ matrix are immediate
consequences of the inequalities involving variance of real and complex
numbers.



\bigskip

\bigskip

\bigskip

\bigskip

\bigskip

\bigskip

\bigskip

\bigskip

\bigskip

\bigskip

\bigskip

\bigskip

\bigskip

\bigskip

\bigskip

\bigskip

\bigskip

\bigskip

\bigskip

\section{Introduction}

\setcounter{equation}{0}\qquad It is useful to have bounds for eigenvalues
and spread in term of the functions of entries of the given matrix. Such
bounds have been studied extensively in literature, see [1-18]. Bounds on
eigenvalues and spread of a matrix $A$ in terms of the traces of $A$ and $%
A^{2}$ are of special interests. These bounds are in fact the immediate
consequences of inequalities for real or complex numbers. In this note we
point out some more such bounds related to the variance of real and complex
numbers.

Let $x_{1},x_{2},...,x_{n}$ denote $n$ real numbers. Their arithmetic mean
is the number%
\begin{equation}
\overline{x}=\frac{1}{n}\sum_{i=1}^{n}x_{i}  \tag{1.1}
\end{equation}%
and the variance is%
\begin{equation}
S_{x}^{2}=\frac{1}{n}\sum_{i=1}^{n}\left( x_{i}-\overline{x}\right) ^{2}. 
\tag{1.2}
\end{equation}%
Samuelson's inequality \cite{12} says that%
\begin{equation}
S_{x}^{2}\geq \frac{1}{n-1}\left( x_{j}-\overline{x}\right) ^{2},  \tag{1.3}
\end{equation}%
for all $j=1,2,...,n$.

\vskip0.0in\noindent The Nagy inequality \cite{11} gives a lower bound for
the variance,%
\begin{equation}
S_{x}^{2}\geq \frac{1}{2n}\max_{i,j}\left( x_{i}-x_{j}\right) ^{2}. 
\tag{1.4}
\end{equation}%
A more general inequality due to Fahmy and Prochan \cite{7} says that for $%
x_{1}\leq x_{2}\leq ...\leq x_{n},$ we have%
\begin{equation}
S_{x}^{2}\geq \frac{l\left( n-k+1\right) }{n\left( n+l-k+1\right) }\left(
x_{k}-x_{l}\right) ^{2},  \tag{1.5}
\end{equation}%
where $1\leq l<k\leq n.$

We here study some extensions of these inequalities for the complex numbers
and discuss their applications. Let $z_{1},z_{2},...,z_{n}$ denote $n$
complex numbers. Their arithmetic mean is the number%
\begin{equation}
\widetilde{z}=\frac{1}{n}\sum_{i=1}^{n}z_{i}.  \tag{1.6}
\end{equation}%
Let%
\begin{equation}
S^{2}=\frac{1}{n}\sum_{i=1}^{n}\left( z_{i}-\widetilde{z}\right) ^{2}, 
\tag{1.7}
\end{equation}%
and%
\begin{equation}
S_{z}^{2}=\frac{1}{n}\sum_{i=1}^{n}\left\vert z_{i}-\widetilde{z}\right\vert
^{2}.  \tag{1.8}
\end{equation}%
For $z_{i}=x_{i}+\imath y_{i}$ we have $S_{z}^{2}=S_{x}^{2}+S_{y}^{2}.$ On
substituting $2x_{i}=z_{i}+\overline{z_{i}}$ in (1.2), we get that%
\begin{equation}
2S_{x}^{2}=S_{z}^{2}+\func{Re}S^{2}\text{ and\ \ }2S_{y}^{2}=S_{z}^{2}-\func{%
Re}S^{2}.  \tag{1.9}
\end{equation}%
We first prove some basic inequalities involving real and complex numbers in
the following lemmas, and use these inequalities to derive several bounds
for the eigenvalues in Section 2.

\bigskip \vskip0.0in\noindent \textbf{Lemma 1.1.} Let $z_{i}=x_{i}+\imath
y_{i},\ i=1,2,...n$. For $x_{1}\leq x_{2}\leq ...\leq x_{n},$%
\begin{equation}
\left\vert x_{k}-x_{l}\right\vert ^{2}\leq \frac{n\left( n+l-k+1\right) }{%
2l\left( n-k+1\right) }\left( S_{z}^{2}+\func{Re}S^{2}\right) ,  \tag{1.10}
\end{equation}%
For $y_{1}\leq y_{2}\leq ...\leq y_{n},$%
\begin{equation}
\left\vert y_{k}-y_{l}\right\vert ^{2}\leq \frac{n\left( n+l-k+1\right) }{%
2l\left( n-k+1\right) }\left( S_{z}^{2}-\func{Re}S^{2}\right) .  \tag{1.11}
\end{equation}%
Also, the inequalities%
\begin{equation}
\left\vert z_{k}-z_{l}\right\vert ^{2}\leq \frac{n\left( n+l-k+1\right) }{%
2l\left( n-k+1\right) }\left( S_{z}^{2}+\left\vert S^{2}\right\vert \right) ,
\tag{1.12}
\end{equation}%
hold for some permutation of complex numbers $z_{i}$, $i=1,2,...,n,$ and $%
1\leq l<k\leq n.$

\vskip0.0in\noindent \textbf{Proof. }For any complex number $z$ there is a
complex number $\alpha $ with $\left\vert \alpha \right\vert =1$ such that $%
\func{Re}(\alpha z)=\left\vert z\right\vert $. Therefore, for $x_{i}=\func{Re%
}(\alpha z_{i})$ we can choose $\alpha $ such that $\left\vert \alpha
\right\vert =1$ and 
\begin{equation}
\left\vert x_{k}-x_{l}\right\vert =\left\vert \func{Re}(\alpha z_{k})-\func{%
Re}(\alpha z_{l})\right\vert =\left\vert z_{k}-z_{l}\right\vert .  \tag{1.13}
\end{equation}%
Without restricting generality, assume that $z_{1},z_{2},...,z_{n}$ is a
permutation of complex numbers $z_{i}$ such that $\func{Re}(\alpha
z_{1})\leq \func{Re}(\alpha z_{2})\leq ....\leq \func{Re}(\alpha z_{n})$.
Put $x_{i}=\func{Re}(\alpha z_{i})$ in (1.2), we find that%
\begin{equation}
2S_{x}^{2}=S_{z}^{2}+\frac{1}{n}\func{Re}(\alpha ^{2}\sum_{i=1}^{n}\left(
z_{i}-\widetilde{z}\right) ^{2}).  \tag{1.14}
\end{equation}%
For $\left\vert \alpha \right\vert =1$, we have 
\begin{equation}
\left\vert \frac{1}{n}\func{Re}(\alpha ^{2}\sum_{i=1}^{n}\left( z_{i}-%
\widetilde{z}\right) ^{2})\right\vert \leq \left\vert S^{2}\right\vert . 
\tag{1.15}
\end{equation}%
It follows from (1.14) and (1.15) that for any complex number $\alpha $ with 
$\left\vert \alpha \right\vert =1$ and $x_{i}=\func{Re}(\alpha z_{i})$, we
have%
\begin{equation}
2S_{x}^{2}\leq S_{z}^{2}+\left\vert S^{2}\right\vert .  \tag{1.16}
\end{equation}%
Substituting (1.13) in (1.5) and use (1.16), we immediately get (1.12).

\noindent The inequalities (1.10) and (1.11) follow immediately on using
(1.9) in (1.5). \ \ $\blacksquare $

\vskip0.0in\noindent \textbf{Lemma 1.2.} Let%
\begin{equation}
\mu _{2}=\sum\limits_{i=1}^{n}p_{i}\left( x_{i}-\mu _{1}^{^{\prime }}\right)
^{2},  \tag{1.17}
\end{equation}%
where $\mu _{1}^{^{\prime }}=\sum\limits_{i=1}^{n}p_{i}x_{i}$ and $p_{i}$
are non-negative real numbers such that $\sum\limits_{i=1}^{n}p_{i}=1$ and $%
p_{j}\neq 1.$ Then%
\begin{equation}
\mu _{2}\geq \frac{p_{j}}{1-p_{j}}\left( \mu _{1}^{^{\prime }}-x_{j}\right)
^{2},  \tag{1.18}
\end{equation}%
for all $j=1,2,\ldots ,n.$

\vskip0.0in\noindent \textbf{Proof. }From (1.17), we have%
\begin{equation}
\mu _{2}=p_{j}\left( x_{j}-\mu _{1}^{^{\prime }}\right) ^{2}+\left(
1-p_{j}\right) \frac{1}{1-p_{j}}\sum_{i=1,i\neq j}^{n}p_{i}\left( x_{i}-\mu
_{1}^{^{\prime }}\right) ^{2}.  \tag{1.19}
\end{equation}%
It follows from the Cauchy-Schwarz inequality that%
\begin{equation}
\frac{1}{1-p_{j}}\sum_{i=1,i\neq j}^{n}p_{i}\left( x_{i}-\mu _{1}^{^{\prime
}}\right) ^{2}\geq \left( \frac{1}{1-p_{j}}\sum_{i=1,i\neq j}^{n}p_{i}\left(
x_{i}-\mu _{1}^{^{\prime }}\right) \right) ^{2}.  \tag{1.20}
\end{equation}%
On the other hand, the sum of all the deviations from the mean is zero,
therefore%
\begin{equation*}
\sum_{i=1}^{n}p_{i}\left( x_{i}-\mu _{1}^{^{\prime }}\right) =0,
\end{equation*}%
and we get that%
\begin{equation}
\sum_{i=1,i\neq j}^{n}p_{i}\left( x_{i}-\mu _{1}^{^{\prime }}\right)
=p_{j}\left( x_{j}-\mu _{1}^{^{\prime }}\right) .  \tag{1.21}
\end{equation}%
Combining (1.20) and (1.21), we find that 
\begin{equation}
\frac{1}{1-p_{j}}\sum_{i=1,i\neq j}^{n}p_{i}\left( x_{i}-\mu _{1}^{^{\prime
}}\right) ^{2}\geq \left( \frac{p_{j}\left( x_{j}-\mu _{1}^{^{\prime
}}\right) }{1-p_{j}}\right) ^{2}.  \tag{1.22}
\end{equation}%
Insert (1.22) in (1.19), a little computation leads to (1.18). \ \ $%
\blacksquare $

\vskip0.0in\noindent \textbf{Lemma 1.3. }For $x_{1}\leq x_{2}\leq \ldots
\leq x_{n},$ the inequality%
\begin{equation}
S_{x}^{2}\geq \frac{k}{n-k}\left( \overline{x}-x_{j}\right) ^{2}  \tag{1.23}
\end{equation}%
holds for $j=k$ and $j=n-k+1$ with $k\leq n-k+1.$

\vskip0.0in\noindent \textbf{Proof. }For $\overline{x}\geq x_{k},$ we have%
\begin{equation}
S_{x}^{2}=\frac{1}{n}\sum_{i=1}^{n}\left( x_{i}-\overline{x}\right) ^{2}\geq 
\frac{k}{n}\left( x_{k}-\overline{x}\right) ^{2}+\frac{1}{n}%
\sum_{i=k+1}^{n}\left( x_{i}-\overline{x}\right) ^{2}.  \tag{1.24}
\end{equation}%
Let%
\begin{equation*}
\overline{y}=\frac{1}{n}\left( kx_{k}+\sum_{i=k+1}^{n}x_{i}\right) ,
\end{equation*}%
and%
\begin{equation*}
S_{y}^{2}=\frac{1}{n}\left( k\left( x_{k}-\overline{y}\right)
^{2}+\sum_{i=k+1}^{n}\left( x_{i}-\overline{y}\right) ^{2}\right) .
\end{equation*}%
It is clear that $\overline{y}\geq \overline{x}$ and 
\begin{equation*}
S_{x}^{2}\geq \frac{k}{n}\left( x_{k}-\overline{x}\right) ^{2}+\frac{1}{n}%
\sum_{i=k+1}^{n}\left( x_{i}-\overline{x}\right) ^{2}\geq S_{y}^{2}.
\end{equation*}%
Using Lemma 1.2, we have%
\begin{equation*}
S_{x}^{2}\geq S_{y}^{2}\geq \frac{k}{n-k}\left( \overline{y}-x_{k}\right)
^{2}\geq \frac{k}{n-k}\left( \overline{x}-x_{k}\right) ^{2},
\end{equation*}%
for $\ k=1,2,...,n-1.$

On using similar argument we find that for $\overline{x}\leq x_{k},$ we have%
\begin{equation*}
S_{x}^{2}\geq \frac{n-k+1}{k-1}\left( x_{k}-\overline{x}\right) ^{2},\text{
\ }
\end{equation*}%
\ $k=2,3,..,n.$ Thus,%
\begin{equation}
S_{x}^{2}\geq \min \left\{ \frac{k}{n-k}\left( \overline{x}-x_{k}\right)
^{2},\ \frac{n-k+1}{k-1}\left( x_{k}-\overline{x}\right) ^{2}\right\} , 
\tag{1.25}
\end{equation}%
for $k=2,3,...,n-1.$ The assertions of lemma now follow from (1.25) and the
fact that $\left( n-k+1\right) \left( n-k\right) \geq k\left( k-1\right) $
if and only if $k\leq n-k+1.$ $\ \ \blacksquare $

\vskip0.0in\noindent \textbf{Lemma 1.4. }Let $z_{i}=x_{i}+\imath y_{i},\
i=1,2,...n$. For $x_{1}\leq x_{2}\leq ...\leq x_{n},$ the inequalities%
\begin{equation}
S_{z}^{2}+\func{Re}S^{2}\geq \frac{2k}{n-k}\left\vert \frac{1}{n}%
\sum_{i=1}^{n}x_{i}-x_{j}\right\vert ^{2},  \tag{1.26}
\end{equation}%
hold for $j=k$ and $j=n-k+1$ with $k\leq n-k+1.$ Likewise, for $y_{1}\leq
y_{2}\leq ...\leq y_{n}$, we have%
\begin{equation}
S_{z}^{2}-\func{Re}S^{2}\geq \frac{2k}{n-k}\left\vert \frac{1}{n}%
\sum_{i=1}^{n}y_{i}-y_{j}\right\vert ^{2}.  \tag{1.27}
\end{equation}

\vskip0.0in\noindent Also,%
\begin{equation}
S_{z}^{2}+\left\vert S^{2}\right\vert \geq \frac{2k}{n-k}\left\vert \frac{1}{%
n}\sum_{i=1}^{n}z_{i}-z_{j}\right\vert .  \tag{1.28}
\end{equation}%
\vskip0.0in\noindent \textbf{Proof. }Use (1.6)-(1.9), (1.16) and Lemma 1.3;
we immediately get the inequalities (1.26)-(1.28).

\vskip0.0in\noindent \textbf{Lemma 1.5.} For $n$ complex numbers $z_{i},$ we
have%
\begin{equation}
\frac{1}{n}\sum_{i=1}^{n}\left\vert z_{i}-\widetilde{z}\right\vert ^{2r}\geq 
\frac{1+\left( n-1\right) ^{2r-1}}{n\left( n-1\right) ^{2r-1}}\left\vert
z_{j}-\widetilde{z}\right\vert ^{2r},  \tag{1.29}
\end{equation}%
where $\widetilde{z}$\ is given in (1.6).

\vskip0.0in\noindent \textbf{Proof. }We write 
\begin{equation}
\frac{1}{n}\sum_{i=1}^{n}\left\vert z_{i}-\widetilde{z}\right\vert ^{2r}=%
\frac{\left\vert z_{j}-\widetilde{z}\right\vert ^{2r}}{n}+\frac{n-1}{n}\frac{%
1}{n-1}\dsum\limits_{i=1,i\neq j}^{n}\left\vert z_{i}-\widetilde{z}%
\right\vert ^{2r}.  \tag{1.30}
\end{equation}%
For $m$ positive real numbers $y_{i}$, $i=1,2,...,m$, 
\begin{equation}
\frac{1}{m}\dsum\limits_{i=1}^{m}y_{i}^{k}\geq \left( \frac{1}{m}%
\dsum\limits_{i=1}^{m}y_{i}\right) ^{k},\text{ }k=1,2,\ldots \text{ .} 
\tag{1.31}
\end{equation}%
Applying (1.31) to $n-1$ positive real numbers $\left\vert z_{i}-\widetilde{z%
}\right\vert ^{2r},i=1,2,\ldots ,n$ and $i\neq j,$ we get%
\begin{equation}
\frac{1}{n-1}\dsum\limits_{i=1,i\neq j}^{n}\left( \left\vert z_{i}-%
\widetilde{z}\right\vert ^{2}\right) ^{r}\geq \left( \frac{1}{n-1}%
\dsum\limits_{i=1,i\neq j}^{n}\left\vert z_{i}-\widetilde{z}\right\vert
^{2}\right) ^{r}.  \tag{1.32}
\end{equation}%
Also, we have%
\begin{equation}
\frac{1}{n-1}\dsum\limits_{i=1,i\neq j}^{n}\left\vert z_{i}-\widetilde{z}%
\right\vert ^{2}\geq \left( \frac{1}{n-1}\dsum\limits_{i=1,i\neq
j}^{n}\left\vert z_{i}-\widetilde{z}\right\vert \right) ^{2}.  \tag{1.33}
\end{equation}%
On the other hand the sum of all the deviations from the mean is zero,
therefore 
\begin{equation*}
\dsum\limits_{i=1}^{n}\left( z_{i}-\widetilde{z}\right) =0,
\end{equation*}%
and we get that%
\begin{equation}
\left\vert \widetilde{z}-z_{j}\right\vert \leq \dsum\limits_{i=1,i\neq
j}^{n}\left\vert z_{i}-\widetilde{z}\right\vert .  \tag{1.34}
\end{equation}%
Combining (1.31)-(1.34), we find that%
\begin{equation}
\frac{1}{n-1}\dsum\limits_{i=1,i\neq j}^{n}\left\vert z_{i}-\widetilde{z}%
\right\vert ^{2r}\geq \left( \frac{\left\vert z_{j}-\widetilde{z}\right\vert 
}{n-1}\right) ^{2r}.  \tag{1.35}
\end{equation}%
Insert (1.35) in (1.30), a little computation leads to (1.29). \ \ $%
\blacksquare $

\section{Bounds on eigenvalues using traces}

\setcounter{equation}{0}\noindent Let $A=\left( a_{ij}\right) $ be an $%
n\times n$ complex matrix with eigenvalues $\lambda _{1},\lambda
_{2},...,\lambda _{n}$. Let $\lambda _{i}=\alpha _{i}+\imath \beta _{i},\
i=1,2,...n$. Then%
\begin{equation}
\widetilde{\lambda }=\frac{1}{n}\sum_{i=1}^{n}\lambda _{i}=\frac{\text{tr}A}{%
n},  \tag{2.1}
\end{equation}%
\begin{equation}
S^{2}=\frac{1}{n}\sum_{i=1}^{n}\left( \lambda _{i}-\widetilde{\lambda }%
\right) ^{2}=\frac{\text{tr}A^{2}}{n}-\left( \frac{\text{tr}A}{n}\right)
^{2},  \tag{2.2}
\end{equation}%
and%
\begin{equation}
S_{\lambda }^{2}=\frac{1}{n}\sum_{i=1}^{n}\left\vert \lambda _{i}-\widetilde{%
\lambda }\right\vert ^{2}=\frac{1}{n}\sum_{i=1}^{n}\left\vert \lambda
_{i}\right\vert ^{2}-\left\vert \frac{\text{tr}A}{n}\right\vert ^{2}, 
\tag{2.3}
\end{equation}%
where tr$A$ denotes the trace of $A$. Likewise, we have%
\begin{equation}
S_{\alpha }^{2}=\frac{1}{n}\sum_{i=1}^{n}\left( \alpha _{i}-\widetilde{%
\alpha }\right) ^{2}=\frac{S_{\lambda }^{2}+\func{Re}S^{2}}{2},  \tag{2.4}
\end{equation}%
and%
\begin{equation}
S_{\beta }^{2}=\frac{1}{n}\sum_{i=1}^{n}\left( \beta _{i}-\widetilde{\beta }%
\right) ^{2}=\frac{S_{\lambda }^{2}-\func{Re}S^{2}}{2}.  \tag{2.5}
\end{equation}

\vskip0.0in\noindent \textbf{Theorem 2.1. }Let $A$ be an $n\times n$ complex
matrix. Then, the disk%
\begin{equation}
\left\vert z-\frac{\text{tr}A}{n}\right\vert \leq \sqrt{\frac{n-k}{2k}\left(
S_{\lambda }^{2}+\left\vert S\right\vert ^{2}\right) },  \tag{2.6}
\end{equation}%
contains at least $n-2k+2$ eigenvalues of $A$ where $k$ is a positive
integer less than or equal to $\frac{n+1}{2}.$ Also, the disks%
\begin{equation}
\left\vert x-\frac{\func{Re}\text{tr}A}{n}\right\vert \leq \sqrt{\frac{n-k}{%
2k}\left( S_{\lambda }^{2}+\func{Re}S^{2}\right) }  \tag{2.7}
\end{equation}%
and%
\begin{equation}
\left\vert y-\frac{\func{Im}\text{tr}A}{n}\right\vert \leq \sqrt{\frac{n-k}{%
2k}\left( S_{\lambda }^{2}-\func{Re}S^{2}\right) },  \tag{2.8}
\end{equation}%
contains real and imaginary parts of eigenvalues, respectively.

\vskip0.0in\noindent \textbf{Proof. }Let\textbf{\ }$x_{i}=\func{Re}\left(
\alpha \lambda _{i}\right) .$ Then, there is a complex number $\alpha $ with 
$\left\vert \alpha \right\vert =1$ such that%
\begin{equation*}
\left\vert x_{i}-\frac{1}{n}\sum_{i=1}^{n}x_{i}\right\vert =\left\vert \func{%
Re}\alpha \left( \lambda _{i}-\frac{\text{tr}A}{n}\right) \right\vert
=\left\vert \lambda _{i}-\frac{\text{tr}A}{n}\right\vert .
\end{equation*}%
It follows from Lemma 1.4 that%
\begin{equation}
\left\vert \lambda _{i}-\frac{\text{tr}A}{n}\right\vert \leq \sqrt{\frac{n-k%
}{2k}\left( S_{\lambda }^{2}+\left\vert S\right\vert ^{2}\right) }, 
\tag{2.9}
\end{equation}%
for $i=k$ and $i=n-k+1$ with $k\leq n-k+1.$ Let $D_{k}$ denote the $k^{th}$
disk in (2.9). It is easy to see that $D_{\frac{n+1}{2}}\subseteq D_{\frac{%
n-1}{2}}\subseteq ...\subseteq D_{2}\subseteq D_{1}$ when $n$ is odd. By
Lemma 1.4, the eigenvalue $\lambda _{\frac{n+1}{2}}$ lies in the disk $D_{%
\frac{n+1}{2}}$ while $\lambda _{\frac{n-1}{2}}$ and $\lambda _{\frac{n+3}{2}%
}$ lie in $D_{\frac{n+1}{2}}\subseteq D_{\frac{n-1}{2}}$. So the disk $D_{%
\frac{n-1}{2}}$ contains at least three eigenvalues. Repeating the above
process, we can easily see that the disk $D_{k}$ contains $n-2k+2$
eigenvalues. Likewise, the inequalities (2.7) and (2.8) follow from (1.26)
and (1.27), respectively. \ \ $\blacksquare $

\vskip0.0in\noindent For a normal matrix, we have%
\begin{equation}
S_{\lambda }^{2}=\frac{\text{tr}AA^{\ast }}{n}-\left\vert \frac{\text{tr}A}{n%
}\right\vert ^{2}.  \tag{2.10}
\end{equation}

Using this in (2.6)-(2.8), we can calculate the corresponding upper bounds
and hence the regions containing eigenvalues. For arbitrary matrices, we can
use the various upper bounds for $\sum_{i=1}^{n}\left\vert \lambda
_{i}\right\vert ^{2},\ $see \cite{17}. For special case $k=1$, we get
Theorem 2.7 of Huang and Wang \cite{8}. It is clear that if rank of matrix
is $m;$ we can replace $n$ by $m$.

\vskip0.0in\noindent \textbf{Theorem 2.2. }Let $A$ be an $n\times n$ complex
matrix with eigenvalues $\lambda _{i}=\lambda _{i}=\alpha _{i}+\imath \beta
_{i},\ i=1,2,...n$. For $\alpha _{1}\leq \alpha _{2}\leq ...\leq \alpha
_{n}, $ the inequalities%
\begin{equation}
\left\vert \alpha _{k}-\alpha _{l}\right\vert ^{2}\leq \frac{n\left(
n+l-k+1\right) }{2l\left( n-k+1\right) }\left( S_{\lambda }^{2}+\func{Re}%
S^{2}\right) ,  \tag{2.11}
\end{equation}%
hold for $1\leq l<k\leq n$, and for $\beta _{1}\leq \beta _{2}\leq ...\leq
\beta _{n}$%
\begin{equation}
\left\vert \beta _{k}-\beta _{l}\right\vert ^{2}\leq \frac{n\left(
n+l-k+1\right) }{2l\left( n-k+1\right) }\left( S_{\lambda }^{2}-\func{Re}%
S^{2}\right) .  \tag{2.12}
\end{equation}%
Also, the inequalities%
\begin{equation}
\left\vert \lambda _{k}-\lambda _{l}\right\vert ^{2}\leq \frac{n\left(
n+l-k+1\right) }{2l\left( n-k+1\right) }\left( S_{\lambda }^{2}+\left\vert
S^{2}\right\vert \right) ,  \tag{2.13}
\end{equation}%
hold for some permutation of complex numbers $z_{i}$, $i=1,2,...,n$ and $%
1\leq l<k\leq n.$

\vskip0.0in\noindent \textbf{Proof. }The assertions of the theorem follow
easily from Lemma 1.1. \ \ $\blacksquare $

\vskip0.1inIt may be noted here that for $l=1$ and $k=n,$ Theorem 2.2
provides bounds for the spread, Spd$(A)=\max_{i,j}\left\vert \lambda
_{i}-\lambda _{j}\right\vert .$

\vskip0.0in\noindent \textbf{Corollary 2.3. }Let $A$ be an $n\times n$
complex matrix with at least two distinct eigenvalues. Let $\lambda
_{k}=\alpha _{k}+\imath \beta _{k},$ $k=1,2,...n$ be any eigenvalue of $A$.
For $\alpha _{1}\leq \alpha _{2}\leq ...\leq \alpha _{n},$ the disk%
\begin{equation}
\left\vert x-\alpha _{l}\right\vert \leq \frac{n}{\sqrt{2\left( n-1\right) }}%
\sqrt{S_{\lambda }^{2}+\func{Re}S^{2}},  \tag{2.14}
\end{equation}%
contains real part of one more eigenvalue. For $\beta _{1}\leq \beta
_{2}\leq ...\leq \beta _{n}$%
\begin{equation}
\left\vert y-\beta _{l}\right\vert \leq \frac{n}{\sqrt{2\left( n-1\right) }}%
\sqrt{S_{\lambda }^{2}-\func{Re}S^{2}},  \tag{2.15}
\end{equation}%
contains imaginary part of one more eigenvalue. Also, the disk%
\begin{equation}
\left\vert z-\lambda _{l}\right\vert \leq \frac{n}{\sqrt{2\left( n-1\right) }%
}\sqrt{S_{\lambda }^{2}+\left\vert S^{2}\right\vert },  \tag{2.16}
\end{equation}%
contains one more eigenvalue of $A.$

\vskip0.0in\noindent \textbf{Proof. }It follows from Theorem 2.2 that for $%
k=l+1,$%
\begin{equation}
2\left\vert \alpha _{l+1}-\alpha _{l}\right\vert ^{2}\leq \frac{n^{2}}{%
2l(n-l)}\left( S_{\lambda }^{2}+\func{Re}S^{2}\right) ,  \tag{2.17}
\end{equation}%
\begin{equation}
2\left\vert \beta _{l+1}-\beta _{l}\right\vert ^{2}\leq \frac{n^{2}}{2l(n-l)}%
\left( S_{\lambda }^{2}-\func{Re}S^{2}\right) ,  \tag{2.18}
\end{equation}%
and%
\begin{equation}
2\left\vert \lambda _{l+1}-\lambda _{l}\right\vert ^{2}\leq \frac{n^{2}}{%
l(n-l)}\left( S_{\lambda }^{2}+\left\vert S^{2}\right\vert \right) , 
\tag{2.19}
\end{equation}%
$l=1,2,...,n-1$. The inequalities (2.14)-(2.16) respectively follow from
(2.17)-(2.19), and the fact that $l(n-l)\geq n-1$ for $l=1,2,...,n-1.$ $\ \
\blacksquare $

\vskip0.0in\noindent The above corollary is useful when one of the
eigenvalues is known. For example, for a singular matrix one eigenvalue is
zero.

\vskip0.0in\noindent \textbf{Theorem 2.4. }Let $A$ be an $n\times n$ complex
matrix. Then all the eigenvalues of $A$ are contained in the disk%
\begin{equation}
\left\vert z-\frac{\text{tr}A}{n}\right\vert \leq \left( \frac{\left(
n-1\right) ^{2r-1}}{1+\left( n-1\right) ^{2r-1}}\sum_{i=1}^{n}\left\vert
\lambda _{i}\left( B\right) \right\vert ^{2r}\right) ^{\frac{1}{2r}}, 
\tag{2.20}
\end{equation}%
where $B=A-\frac{\text{tr}A}{n}I$ and $r=1,2,....$

\vskip0.0in\noindent \textbf{Proof.} Let $\lambda _{i}$ be eigenvalues of $A$%
, then the eigenvalues of $B^{r}$ are $\left( \lambda _{i}-\frac{\text{tr}A}{%
n}\right) ^{r},\ i=1,2,...n$. The assertions of the theorem now follow on
using Lemma 1.5. \ \ $\blacksquare $

Note that%
\begin{eqnarray*}
\sum_{i=1}^{n}\left\vert \lambda _{i}\left( B\right) \right\vert ^{2r}
&=&\sum_{i=1}^{n}\left\vert \lambda _{i}\left( B\right) ^{r}\right\vert ^{2}
\\
&\leq &\sqrt{\left( \left\Vert B^{r}\right\Vert _{2}^{2}-\frac{1}{n}%
\left\vert \text{tr}B^{r}\right\vert ^{2}\right) ^{2}-\frac{1}{2}\left\Vert
B^{r}\left( B^{r}\right) ^{\ast }-\left( B^{r}\right) ^{\ast
}B^{r}\right\Vert _{2}^{2}}+\frac{1}{n}\left\vert \text{tr}B^{r}\right\vert
^{2}\text{.}
\end{eqnarray*}%
Using this and other relation in (2.20), we can obtain various bounds for
eigenvalues.

\vskip0.0in\noindent \textbf{Theorem 2.5. }Let $A$ be an $n\times n$ complex
normal matrix. Then, one eigenvalue of $A$ lies on or outside the circle%
\begin{equation}
\left\vert z-\frac{\text{tr}A}{n}\right\vert =\left( \text{tr}B^{r}\left(
B^{r}\right) ^{\ast }\right) ^{\frac{1}{2r}}.  \tag{2.21}
\end{equation}

\vskip0.0in\noindent \textbf{Proof.} The proof of theorem follows from the
fact that%
\begin{equation*}
\frac{\text{tr}B^{r}\left( B^{r}\right) ^{\ast }}{n}=\frac{1}{n}%
\sum_{i=1}^{n}\left\vert \lambda _{i}-\frac{\text{tr}A}{n}\right\vert
^{2r}\leq \left\vert \lambda _{j}-\frac{\text{tr}A}{n}\right\vert ^{2r}.\ \
\ \blacksquare
\end{equation*}

\vskip0.0inLet the eigenvalues of a complex $n\times n$ matrix $A$ are all
real, as in case of a Hermitian matrix. Wolkowicz and Styan \cite{18} have
shown that%
\begin{equation*}
\lambda _{\max }\geq \frac{\text{tr}A}{n}+\sqrt{\frac{\text{tr}B^{2}}{%
n\left( n-1\right) }}
\end{equation*}%
and%
\begin{equation*}
\lambda _{\min }\leq \frac{\text{tr}A}{n}-\sqrt{\frac{\text{tr}B^{2}}{%
n\left( n-1\right) }}.
\end{equation*}%
We prove extensions of these inequalities in the following theorem.

\vskip0.0in\noindent \textbf{Theorem 2.6.} Let $A$ be an $n\times n$ complex
matrix with real eigenvalues $\lambda _{i},i=1,2,...,n.$ Then%
\begin{equation}
\lambda _{\max }\geq \frac{\text{tr}A}{n}+\frac{\text{tr}B^{2}}{n}\left( 
\frac{1+\left( n-1\right) ^{2r-1}}{\left( n-1\right) ^{2r-1}}\frac{1}{\text{%
tr}B^{2r}}\right) ^{\frac{1}{2r}}  \tag{2.22}
\end{equation}%
and%
\begin{equation}
\lambda _{\min }\leq \frac{\text{tr}A}{n}-\frac{\text{tr}B^{2}}{n}\left( 
\frac{1+\left( n-1\right) ^{2r-1}}{\left( n-1\right) ^{2r-1}}\frac{1}{\text{%
tr}B^{2r}}\right) ^{\frac{1}{2r}},  \tag{2.23}
\end{equation}%
where $B=A-\frac{\text{tr}A}{n}I.$

\vskip0.0in\noindent \textbf{Proof. }We have \cite{3}%
\begin{equation*}
\frac{1}{n}\text{tr}B^{2}\leq \left( \lambda _{\max }-\frac{\text{tr}A}{n}%
\right) \left( \frac{\text{tr}A}{n}-\lambda _{\min }\right) .
\end{equation*}%
Therefore,%
\begin{equation}
\lambda _{\max }\geq \frac{\text{tr}A}{n}+\frac{1}{n}\frac{\text{tr}B^{2}}{%
\frac{\text{tr}A}{n}-\lambda _{\min }}.  \tag{2.24}
\end{equation}%
It follows from Theorem 2.4 that%
\begin{equation}
\lambda _{\max }\leq \frac{\text{tr}A}{n}+\left( \frac{\left( n-1\right)
^{2r-1}}{1+\left( n-1\right) ^{2r-1}}\text{tr}B^{2r}\right) ^{\frac{1}{2r}}.
\tag{2.25}
\end{equation}%
Combining (2.24) and (2.25), we get (2.23). The inequality (2.22) follows on
using similar arguments. \ \ $\blacksquare $

\vskip0.0in\noindent \textbf{Example 1. }Let%
\begin{equation*}
A=\left[ 
\begin{array}{cccc}
4 & 0 & 2 & 3 \\ 
0 & 5 & 0 & 1 \\ 
2 & 0 & 6 & 0 \\ 
3 & 1 & 0 & 7%
\end{array}%
\right] .
\end{equation*}%
The estimates of Wolkowicz and Styan \cite{18} give $\lambda _{\max }\geq
7.1583$ and $\lambda _{\min }\leq 3.841$ while our estimates (2.22) and
(2.23) for $r=2$ give $\lambda _{\max }\geq 7.2586$ and $\lambda _{\min
}\leq 3.7414,$ respectively.

\vskip0.1in\noindent \textbf{Acknowledgements}. The authors are grateful to
Prof. Rajendra Bhatia for the useful discussions and suggestions. The first
two authors thank I.S.I. Delhi for a visit in January 2014 when this work
had begun.

\setcounter{equation}{0}

\end{document}